\documentclass{amsart}
\usepackage[dvips]{graphicx}
\usepackage{amscd}
\usepackage{amsmath}
\usepackage{amsxtra}
\usepackage{amsfonts}
\usepackage{amssymb}
\newtheorem{theorem}{Theorem}[section]
\newtheorem{corollary}[theorem]{Corollary}

\newtheorem{proposition}[theorem]{Proposition}
\theoremstyle{definition}
\newtheorem{definition}[theorem]{Definition}
\newtheorem{remark}[theorem]{Remark}

\theoremstyle{remark}

\renewcommand{\theclaim}{\textup{\theclaim}}

\newtheorem*{acknowledgements}{Acknowledgements}

\numberwithin{equation}{section}

\input cyracc.def

\newcommand{\triunghi}[7]{
\begin{table}[ht]
\begin{center}
\begin{tabular}{ccccccccc}

$ $&$ $&$ $&$ $&$#1$&$ $&$ $&$ $&$ $\\
$ $&$ $&$ $& $\nearrow $ &$ $&$\nwarrow $&$ $&$ $&$ $\\
$ $&$ $&$#4 $&$ $&$ $&$ $&$#5 $&$ $&$ $\\
$ $&$ $&$ $&$\nwarrow $&$ $&$\nearrow$&$ $&$ $&$ $\\
$ $&$\swarrow $&$ $&$ $&$#7$&$ $&$ $&$\searrow $&$ $\\
$ $&$ $&$ $&$ $&$\downarrow$&$ $&$ $&$ $&$ $\\

$#2 $&$ $&$\longleftarrow$&$ $&$#6$&$ $&$\longrightarrow$&$ $&$#3 $\\
\end{tabular}
\end{center}
\end{table}
}

\newcommand{\diagrama}[6]{
\begin{table}[htdp]
\begin{center}
\begin{tabular}{cccc}
$ $&$ $&$#4$&$ $\\
$ $&$#1$&$\longrightarrow$&$#2$\\
$#5$&$\downarrow$&$\nearrow$&$ $\\
$ $&$#3$&$#6$&$ $
\end{tabular}
\end{center}
\end{table}
}

\begin{document}
\title[Triangle joins of SP1 equivalence relations and their cost]{Triangle joins of SP1 equivalence relations and their cost}
\author{Teodor \c Stefan B\^ildea}
\address{Department of Mathematics\\
The University of Iowa\\
14 MacLean Hall\\
Iowa City, IA 52242-1419\\
U.S.A.\\} \email{sbildea@math.uiowa.edu}


\begin{abstract}
Free joins of measure preserving equivalence relations, with or without amalgamation were introduced by Gaboriau, who studied their cost and obtained a formula for the cost of free products of groups. We introduce the notion of triangle join of measure preserving equivalence relations and study their cost. As a consequence we compute the cost of Thompson's group $G_{2,1}$, viewed as the minimal amalgam of a triangle of finite groups.
\end{abstract}

\maketitle \tableofcontents

\section{\label{Intro}Introduction}

In his seminal paper \cite{GAB}, D. Gaboriau introduced a numerical quantity associated  to standard measure preserving (SP1) equivalence relations on standard non-atomic Borel probability spaces. He then extended this notion to at most countable groups, via their free actions on standard non-atomic Borel probability spaces. In his paper, Gaboriau studies the cost of several constructions with relation : free amalgamated joins (or free amalgamated products), HNN extensions. He gives formulas for costs and uses this formulas to evaluate the cost of the corresponding group constructions. This paper is aimed at yet another construction involving SP1 equivalence relations, namely the triangle join. This is a generalized version of the free join and the typical example is given by the relation induced by the free action of a realizable triangle amalgam of groups. As an example, we estimate the cost of Thompson's group $G_{2,1}$, realized by Brown in \cite{Brown} as a minimal realizable triangle amalgam of finite groups.
Among other important results, D. Gaboriau proved that countable amenable groups all have cost $1$. An open question in the theory of groups concerns the amenability of Thompson's groups $F_{2,1}\subset  G_{2,1}$. For an introduction to this groups see for example \cite{Brown0},\cite{GH} or  \cite{scott}. 
In this article we use $G_{2,1}$ for the  original group $G$ introduced by Thompson and described by K.S. Brown in \cite{Brown} in sections 2 and 3. 
\par
D. Gaboriau showed that the cost of $F_{2,1}$=1. In this article we show that the cost of $G_{2,1}$ is also 1.
\par The article is organized as follows : in section \ref{free prod groups} we introduce the generalized free product and results on existence of amalgams of groups. Next (section \ref{eq rel}) we introduce standard measure-preserving  equivalence relations, graphings for such relations and the cost. In section \ref{free join} 
the notion of triangle joins of equivalence relations is introduced, starting from the relation generated by a free action of a realizable triangle amalgam of countable groups. We reduce the study of such joins to that of minimal joins and then give a sufficient condition for a triangle join to be a free amalgamated join. The typical examples are the relations generated by free actions of minimal triangle amalgams. In the last section we study the cost of a triangle amalgam of groups. As an application we present  an estimation of the cost of Thompson's Group $G_{2,1}$, represented as a minimal triangle amalgam of finite groups:

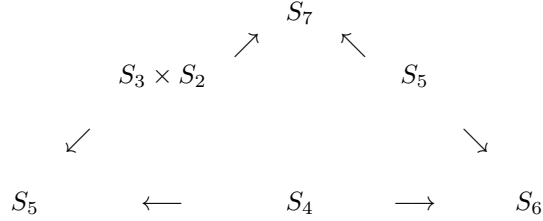
\begin{table}[ht]
\begin{center}
\begin{tabular}{ccccccccc}

$ $&$ $&$ $&$ $&$S_{7}$&$ $&$ $&$ $&$ $\\
$ $&$ $&$ $& $\nearrow $ &$ $&$\nwarrow $&$ $&$ $&$ $\\
$ $&$ $&$S_{3}\times S_2 $&$ $&$ $&$ $&$S_5 $&$ $&$ $\\
$ $&$ $&$ $&$ $&$ $&$ $&$ $&$ $&$ $\\
$ $&$\swarrow $&$ $&$ $&$$&$ $&$ $&$\searrow $&$ $\\
$ $&$ $&$ $&$ $&$ $&$ $&$ $&$ $&$ $\\
$S_5 $&$ $&$\longleftarrow$&$ $&$S_{4}$&$ $&$\longrightarrow$&$ $&$S_{6} $\\
\end{tabular}
\end{center}
\end{table}

K.S. Brown used this presentation of $G_{2,1}$ in \cite{Brown} to settle a couple of questions. One of them is related to the Gersten-Stallings condition (\cite{Sta}) for the realizability of such an amalgam. They defined an angle between the subgroups $H_{12},H_{13}\subset G_1$ as follows. Consider the unique morphism $\phi:H_{12}\ast_K H_{13}\to G_1$. Let $2n=\inf \{|w| | w\in ker\phi,w\ne e\} $ and define $\theta_1=\frac{\pi}{n}$ to be the angle between the two subgroups. Consider the 3 angles involved in the triangle above, and call them $\theta_{1,2,3}$. The result of Gersten  and Stallings states that , if $\theta_1+\theta_2+\theta_3\le \pi$ then the considered triangle amalgam is realizable (i.e  the vertex groups embed in the direct limit of the diagram) . The triangle above is realizable without satisfying this condition: the angles, as computed by Brown, are $\pi/3$ at $S_5$, $\pi/2$ at  $S_6$ and $\pi/3$ at $S_7$. We will use this presentation of $G_{2,1}$ to show that it has cost $1$.
\section{\label{free prod groups}Generalized Amalgams of Groups over Subgroups}
\par

In what follows we present the definition of the generalized amalgam of groups over a family of subgroups. The interested reader is referred to \cite{HN1},\cite{HN2},\cite{BHN} for the initial papers on the subject. For more recent results on group amalgams see \cite{STA} and \cite{Tits}.
\begin{definition}
Let $(G_i)_{i\in I}$ be a family of groups and for each $i\in I$ let $(H_{ij})_{j\in I,j\ne i }$ be a family of subgroups of $G_i$ as above. The generalized group free product of the family $(G_i)_{i\in I}$ with amalgamated subgroups $(H_{ij})_{i,j\in I,i\ne j}$ is the unique ( up to isomorphism ) group $G$ together with homomorphisms $\psi_i:G_i\to G$ satisfying the properties:
\newline (i) the following diagram is commutative:

\diagrama{H_{ij}}{G}{H_{ji}}{\psi_i}{\varphi_{ij}}{\psi_j}
for every $i,j\in I,i\not=j;$
\newline (ii) given any group $K$ and homomorphisms $\phi_i:G_i\to K$ with commutative diagrams

\diagrama{H_{ij}}{K}{H_{ji}}{\phi_i}{\varphi_{ij}}{\phi_j}
for every $i,j\in I,i\not=j$, there exists a unique homomorphism $\Phi:G\to K$ so that the following diagram commutes:

\diagrama{G_i}{G}{K}{\psi_i}{\phi_i}{\Phi}

We will use the notation \[ G=\ast_{i\in I}(G_i,(H_{ij})_{j\ne i}) .\]
When $G_i=\vee_{j\in I,j\ne i}H_{ij}$, we call $G$ the minimal generalized free product or the minimal amalgam of the family $(A_i)_{i\in I}$.
\end{definition}

\begin{remark}
We will also use the following alternative description of the family of groups and subgroups that we will amalgamate : let $G_i$ be groups for $i\in I$ and for each $i,j\in I,i\ne j$ let $H_{ij}(=H_{ji})$ be groups and $\varphi_{ij}:H_{ij}\to G_i$ be injective group homomorphisms. The groups $H_{ij}$ are now abstract groups that have isomorphic copies inside the bigger groups $G_i$ and $G_j$, and this copies will be identified. If we let 
$ K_{ij}:=\varphi_{ij}(H_{ij}) \stackrel{\varphi_{ji}\circ \varphi_{ij}^{-1}}{\simeq} K_{ji}:=\varphi_{ji}(H_{ij})$,
then the family of groups and subgroups $(G_i)_{i\in I}$ , $(K_{ij})_{j\in I,j\ne i}$ is as in the previous definition. By the generalized free product of the family $(G_i)_{i\in I}$ with amalgamated subgroups $(H_{ij})_{j\in I ,j\ne i}$ we will mean $\ast_{i\in I}(G_i,(K_{ij})_{j\in I,j\ne i})$.  Of course we can go both ways with the description of the families and we will use both notation alternatively in this paper.
\end{remark}

\par The first criterion for the existence of the isomorphic copies of the groups $G_i$ inside $G$ refers to the reduced amalgam. For $i\in I$ let $H_i$ be the subgroup of $G_i$ generated by the family $(H_{ij})_{j\in I,j\ne i}$, i.e. $H_i:=\vee_{j\in I,j\ne i}H_{ij} \subseteq  G_i$. The following reduction theorem is due to Hanna Neumann (\cite{HN1}) and has also a nice presentation in \cite{BHN}, section 15. We present the proof for consistency.

\begin{theorem}
With the above notations, the generalized free product G of the family $(G_i)_{i\in I}$ with amalgamated subgroups $(H_{ij})_{i,j\in I,j\ne }$ is realizable if and only if the generalized free product H of the family $(H_i)_{i\in I}$ with amalgamated subgroups $(H_{ij})_{i,j\in I,i\ne })$ is realizable.
\end{theorem} 

{\it Proof:   }
It is clear that if $G$ is realizable, then each $H_i$ will have isomorphic copies as subgroups of the isomorphic copies of $G_i$, hence $H$ is obviously realizable. Suppose now that $H$ is realizable. For each $i\in I$ we first take the free product  $G_i\ast_{H_i}H$.  Next consider 
\[ G=\ast_{H,i\in I}(G_i\ast_{H_i}H).\]
Then each $G_i$ has an isomorphic copy in $G$. Furthermore, since $H$ is the amalgamated subgroup,  and it is the realizable minimal amalgam, we have, using the same notation for the copies of $G_i$ inside $G$ :  $G_i\cap G_j=G_i\cap G_j\cap H=H_i\cap H_j=H_{ij}=H_{ji}\subset G$. It follows that the $G_i's$ embed in $\vee_{i\in I}G_i\subset G$ such that the resulting diagrams concerning the subgroups commute, hence $\ast_{i\in I}(G_i,(H_{ij})_{j\ne i})$ is realizable. 
\newline
Notice that actually $G=\vee_{i\in I}G_i$. Indeed, a moments thought suffices to see that an element $x\in G$ has the form
\[ x=w_1(\{H_i|i\in I\})g_1w_2(\{H_i|i\in I\})g_2...w_n(\{H_i|i\in I\})g_n,\]
with $n\geq 1$ and where $g_k\in G_{i_k},i_1\ne i_2\ne...\ne i_n$ and $w_k(\{H_i | i\in I\})$ are words over the "alphabet" $\cup_{i\in I}H_i$ so that the letters adjacent to an element of $G_{i_k}$ do not belong to $\cup_{j\in I,j\ne i_k}H_{i_k,j}$. But $w_k(\{H_i | i\in I\})\in \vee_{i\in I}G_i$, hence $x\in \vee_{i\in I}G_i$.
\rule{5 pt}{5 pt}

\begin{remark}
We will be interested in this paper in a particular kind of family of groups, that can be represented on a triangle diagram. 
\triunghi{G_1}{G_2}{G_3}{H_{12}}{H_{13}}{H_{23}}{H_{123}}
All arrows are injective group homomorphisms. The addition of the group $H_{123}$ is just a necessary condition for the generalized free product to be realizable. Indeed, if $\ast_{i\in I}(G_i,(H_{ij})_{i\ne j})$, $I=\{ 1,2,3\}$, is realizable, then we may assume that the arrows are inclusions and that $G_i\cap G_j=H_{ij}$. Then $G_1\cap G_2\cap G_3=H_{12}\cap H_{13}=H_{12}\cap H_{23}=H_{13}\cap H_{23}=H_{123}$. Therefore we will assume the triangles to be fillable, i.e. the intersections of pairs of edge groups are isomorphic under the given family of maps; for example 
\[ H_{123}\simeq \varphi_{12}(H_{12})\cap \varphi_{13}(H_{13})(\subset G_1)\stackrel{\varphi_{31}\circ \varphi_{13}^{-1}}{\simeq}\varphi_{32}(H_{23})\cap \varphi_{31}(H_{13})\subset G_3.\]
\end{remark}

\begin{definition} Suppose we are given a triangle of groups and injective group morphisms as the one above, with the additional hypothesis that the vertex groups are generated by the pairs of images of the edge groups - we will call such a triangle a minimal triangle. If the corresponding generalized free product with amalgamations exists, we say the triangle is realizable.
\end{definition}
\begin{remark} The amalgamated free product of groups contains isomorphic copies of each of its  "factors". For the triangle amalgam described above this is not true in general, and the following was proved to be a collapsing triangle of groups:
\[ G_1=<a,b|a^{-1}ba=b^2> \]
\[ G_2=<a,c|c^{-1}ac=a^2> \]
\[ G_3=<b,c|b^{-1}cb=c^2> \]
\[ H_{ij}=G_i\cap G_j\cong \mathbb{Z},1\le i<j\le3.\]
We reproduce here  the proof in \cite{BHN}, chapter V, section 23. The triangle amalgam has generators $a,b,c$ satisfying all three relations above. The first relation gives $b^{-1}ab=ab^{-1}$, hence $b^{-i}ab^i=ab^{-i},i\ge 1$. By symmetry $c^{-i}bc^i=bc^{-i}$. Conjugating the second relation with $b$ gives on the left
\[ b^{-1}c^{-1}acb=c^{-2}ab^{-1}c^2=a^4c^2b^{-1},\]
and on the right
\[b^{-1}a^2b=(ab^{-1})^2=a^2b^{-3}.\]
Hence $c^2=a^{-2}b^{-2}$. But 
\[ bc^2=c^2bc^{-2}=a^{-2}b^{-2}\cdot b\cdot b^2a^2=b^4;\]
thus $c^2=b^3$, and $b,c^2$ commute. But then $c=c^2$ so $c=1$ and also $a=b=1$.
\end{remark}
The last result of this section addresses the realizability of triangles of groups as the one above. 
\begin{theorem}\label{retri}
Given a minimal triangle of groups and injective groups homomorphisms as the previous one, each of the following is a sufficient condition for the triangle to be realizable:
\newline
(i) One of the vertex groups, $G_3$ say, is the free product of $H_{23},H_{13}$ with amalgamation over $H_{123}$, and with inclusion mappings $\varphi_{32}.\varphi_{13}$. In this case we have $\ast_{i=1}^3(G_i,(H_{ij})_{j\ne i})=G_1\ast_{H_{12}}G_2$. 
\newline
(ii) Two of the the groups $G_i$, $G_1$ and $G_2$ say, have the property that every element of $H_{12}$ commutes with every element of $H_{13}$ and of $H_{23}\quad (i=1,2)$. In this case $\ast_{i=1}^{3}(G_i,(H_{ij})_{j\ne i})$ is a quotient of $G_1\ast_{H_{12}}G_2$. 
\end{theorem}

This theorem and its proof can be found in \cite{HN1}, section 9. The proof is not difficult, the idea is to look at the way $H_{13}\ast_{H_{123}}H_{23}$ embeds into $G_1\ast_{H_{12}}G_2$. This is the way to prove similar results for triangles of equivalence relations.

\section{\label{eq rel}SP1 Equivalence Relations}

For a reference to what follows see \cite{GAB}.
\begin{definition} Let $(X,{\mathcal B},\mu)$ be a standard Borel space with  a non-atomic measure $\mu$. We call a equivalence relation $R\subset X\times X$ an {\bf SP equivalence relation} if it satisfies the following properties:
\begin{itemize}
\item ({\bf S}tandard) For each $x\in X$, the equivalence class $R[x]$ of x is countable and $R\subset X\times X$ is a Borel subset of the product $\sigma$-algebra.
\item ({\bf P}reserving) $R$ preserves the measure $\mu$, i.e. all Borel automorphisms of $X$ with the graph inside $R$ preserve $\mu$.
\end{itemize}
If the measure $\mu$ is a probability measure (i.e. $(X,{\mathcal B},\mu)$ is isomorphic to $[0,1]$ with the induced Borel structure and Lebesgue measure), we call an SP equivalence relation $R$ an {\bf SP1 equivalence relation}.
\end{definition}
\begin{definition}
The SP relation $S$ is a subrelation of the SP relation $R$ if for a.e. $x\in X$ we have $S[x]\subset R[x]$. If $R_1,R_2$ are two SP relations on $(X,\mu)$, we denote by $R_1\vee R_2$ the smallest SP relation with the property that for all $x\in X$ the orbit of $x$ contains $R_1[x]\cup R_2[x]$. We call $R_1\vee R_2$ the SP relation generated by $R_1$ and $R_2$.
\end{definition}
\par The main examples for SP1 equivalence relations come from free actions of groups on a standard Borel space endowed with a non-atomic measure. As the free amalgamated product of groups provides the prototype for the free amalgamated  join, the triangle amalgam provides the blueprint for the triangle join of SP1 equivalence relations.

\section{\label{free join}Free Join and Triangle Join of SP1 Equivalence Relations}

Let $R$ be a SP1 relation on $(X,\mu)$ and let $R_1,R_2$ be two subrelations that generate $R$. Suppose $R_3$ is a common subrelation of $R_1$ and $R_2$.
\begin{definition} A sequence $x_1,x_2,...,x_n$ of points in $X$ is called a {\bf reduced sequence} if 
\begin{itemize}
\item each $(x_i,x_{i+1})$ belongs to one of the factors $R_1,R_2$,
\item two succesive paisr $(x_i,x_{i+1})$ belong to distinct factors,
\item if $n>2$ no $(x_i,x_{i+1})$ belongs to $R_3$,
\item if $n=2$ then $x_i\ne x_{i+1}$.
\end{itemize}
\end{definition}

\begin{definition} $R$ is called the free amalgamated join of $R_1$ and $R_2$ over $R_3$, denoted
\[ R=R_1\ast_{R_3}R_2,\]
if every reduced sequence in $X$ has distinct endpoints  up to a set of measure zero).
\end{definition}
\begin{remark} As a typical exmaple consider  the SP1 relation induced by the free action of a free amalgamated product of groups $G=G_1\ast_{G_3}G_2$ on some standard,non-atomic probability space $(X,\mu)$ as above.
\end{remark}

\par Let $G$ be a realizable triangle amalgam of countable groups and let it act freely on a non-atomic measure space. In order to see the defining property of the induced SP1 equivalence relation, we need to exploit the fact that relations in $G$ involve only elements of at most two edge groups and the corresponding vertex group plus identifications of the copies of the edge groups. In other words, the redundancies appear only at the vertices. 

\par Consider the following triangle of SP1 equivalence relations.

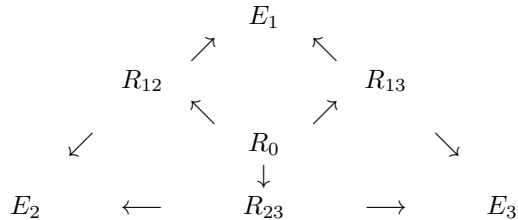
\triunghi{E_1}{E_2}{E_3}{R_{12}}{R_{13}}{R_{23}}{R_0}
Here $R_{ij}=E_i\bigcap E_j$ for $1\le i<j \le 3$.
Let $R_1=R_{12}\vee R_{13}, R_2=R_{12}\vee R_{23}$ and finally $R_3=R_{13}\vee R_{23}$. The following subsets will play a significant role in what follows:
\[ E^{\circ}_i:=E_i\setminus R_i,\mbox{ for  }i=1,2,3;\]
\[  R^{\circ}_{ij}=R_{ij}\setminus R_0\mbox{ for  }1\le i<j\le 3;\]
\[ R=\vee_{i=1}^{3} R_i.\]
The goal is to describe the redundacies in $E_1\vee E_2\vee E_3$. Let us record the following simple facts about the above subsets:
\begin{proposition}\label{absorption}
(i)  {\bf Symmetry:} If $(x,y)\in E^{\circ}_1$, then $(y,x)\in E^{\circ}_1$. If $(x,y)\in R^{\circ}_{12}$, then $(y,x)\in R^{\circ}_{12}$.\\
(ii)  {\bf Absorption 1:} If  $(x,y)\in E^{\circ}_1$ and $(y,z)\in R^{\circ}_{12}$, then $(x,z)\in E^{\circ}_1$;\\
(iii)  {\bf Absorption 2:} If  $(x,y)\in R^{\circ}_{12}$ and $(y,z)\in R_{0}$, then $(x,z)\in R^{\circ}_{12}$.\\ Similar statements hold for $E_2,E_3$ and $R_{13},R_{23}$.
\end{proposition}
{\it Proof:  } Follows from the above definitions.
\rule{5 pt}{5 pt} 
\begin{definition}
(i) A sequence $x_1,x_2,...,x_n\in X$ is called triangle-reduced if for $k\in \{0,1,...,n-1\}$:
\begin{itemize}
\item $(x_k,x_{k+1})\in \sqcup_{i=1}^{3}E_i^{\circ}\sqcup  \sqcup_{1\le i<j\le3}R_{ij}^{\circ}$.
\item $(x_k,x_{k+1})\in E^{\circ}_l \Rightarrow (x_{k-1},x_k),(x_{k+1},x_{k+2})\in E_i^{\circ}\sqcup E_j^{\circ}\sqcup R^{\circ}_{ij}$ where $1\le i<j\le3$ so that $l\ne i,l\ne j$, $l=1,2,3$. 
\item $(x_k,x_{k+1})\in R^{\circ}_{12} \Rightarrow (x_{k-1},x_k),(x_{k+1},x_{k+2})\in E_3^{\circ}\sqcup R_{13}^{\circ}\sqcup R^{\circ}_{23}$. Similar statements for $R_{13}$, $R_{23}$. At endpoints only the left, respectively right neighboring pair is to be considered.
\end{itemize}
(ii) A triangle reduced sequence is called proper if it is of the following types:
\begin{itemize}
\item {\bf type I:} the relations describing the pairs of consecutive points involve (at least) one of the pairs  \[ E_1^{\circ},R^{\circ}_{23} \mbox{ or }E_2^{\circ},R^{\circ}_{13}\mbox{ or }E_3^{\circ}, R^{\circ}_{12}\mbox{ or } \] \[E_1^{\circ},E^{\circ}_{2}\mbox{ or }E_2^{\circ}, E^{\circ}_{3}\mbox{ or }E_3^{\circ}, E^{\circ}_{1}.\]
\item {\bf type II:} the relations describing the pairs of consecutive points involve all  of $R^{\circ}_{12},R^{\circ}_{13},R^{\circ}_{23}$ , but none of $E^{\circ}_1,E^{\circ}_2,E^{\circ}_3$.
\end{itemize}
\end{definition}

\begin{definition}
With the above notations, we say that $E_1\vee E_2\vee E_3$ is a triangle join over $R_{12},R_{23},R_{13}$ if  the proper triangle reduced sequences of type I all have distinct endpoints and the proper triangle reduced loops of type II are concatenations of $R_i$-loops, $i=1,2,3$. 
\end{definition}
\begin{remark} We can also define a "a.e." version of the triangle join, taking into account the probability structure on $(X,\mu)$. Since the triangle join makes sense also on a "bare" set $X$, we will skip the proofs for the "a.e." versions of the next results, noting that they still hold in that framework .
\end{remark}


\begin{remark}
Two extreme cases need to be singled out:
\begin{itemize} 
\item {\bf free amalgamated join:} $R_{12}=R_{23}=R_{13}=R_0=R \not\subseteq E_i$ for $i=1,2,3$. This case corresponds to the free join introduced by Gaboriau (see \cite{GAB} and \cite{KCH}). Indeed, the only triangle reduced sequences are proper of  type I and they correspond to the reduced sequences used by Gaboriau to define the free join.
\item {\bf minimal triangle join:} $E_i=R_i,1\le i \le3$. There are no proper triangle reduced sequences of type I. There is no "extra" information at the vertices, coming from some $E^{\circ}_i$.
\end{itemize}
We will see in what follows what kind of redundancies we expect in a triangle amalgam and how the above cases totally describe all possible situations: up to the minimal join, a triangle join is a free amalgamated join. 
\end{remark}
\begin{proposition}
Suppose $E_1\vee E_2\vee E_3$ is a triangle join over $R_{12},R_{23},R_{13}$. Up to a set of measure zero, between any two distinct points in $X$ there is at most one proper triangle reduced sequence of type one involving only:
\begin{itemize}
\item[(i)] $E^{\circ}_1,E^{\circ}_2$ or $E^{\circ}_2,E^{\circ}_3$ or $E^{\circ}_1,E^{\circ}_3$ or all three of them;
\item[(ii)] $E^{\circ}_1,R^{\circ}_{23}$ or $E^{\circ}_2,R^{\circ}_{13}$ or $E^{\circ}_3,R^{\circ}_{12}$.
\end{itemize}
\end{proposition}

\begin{remark} Before proving the above statement, note that in a triangle join, a proper triangle reduced sequence of type II between $x$ and $y$ is not necessarily unique. This happens because there is no information about the loops at the vertices, or the $R_i$-loops.
\end{remark}
{\it Proof: } For example, in situation (ii), the existence of two distinct sequences with the asserted properties would lead to the existence of a proper triangle reduced loop of type I, which does not exist in a triangle join.
\rule{5 pt}{5 pt}

The next theorem asserts that, in general, a triangle join is a free join up to the minimal triangle join. This shifts the focus to the study of the minimal triangle join. A similar situation is encountered in the study of triangle amalgams of groups or of C*-algebras (see \cite{TSB}).
\begin{theorem}(Reduction to minimal triangle join)
Suppose $E_1\vee E_2\vee E_3$ is a triangle join over $R_{12},R_{23},R_{13}$. Let $R=R_1\vee R_2\vee R_3=R_{12}\vee R_{23}\vee R_{13}$  be the minimal triangle join ($R_i=R_{ij}\vee R_{ik}$). 
Then
\[ E=\ast_{R,i=1}^3(E_i\ast_{R_i}R) \] where $\ast$ is the free (amalgamated) join as defined by Gaboriau in \cite{GAB}.
\end{theorem}

\vspace{0.3 cm}
{\it Proof: }For any equivalence  relation $E$ we will use the notation :
\[ (x,y)\in E\Longleftrightarrow: x\stackrel{\tiny E}{\sim}y .\]
First let us check that 
\begin{equation}\label{l1}E_1\vee R=E_1\ast_{R_1}R.
\end{equation} 
Similar proofs work for $i=2,3$. Let 
\[ x=x_0\stackrel{\tiny E_1}{\sim}x_1\stackrel{\tiny R}{\sim}x_2...\stackrel{\tiny E_1}{\sim}x_{2n-1}\stackrel{\tiny R}{\sim}x_{2n}=x\]
be a sequence in $X$. To prove (\ref{l1}) we need to show that there is some $i,0\le i<2n$ such that $(x_i,x_{i+1})\in R_1$. Assume there is no such $i$. Then $(x_{2k},x_{2k+1})\in E^{\circ}_1$ and $(x_{2k+1},x_{2k+2})\in R\setminus R_1=R_{12}\vee R_{13}\vee R_{23}\setminus R_{12}\vee R_{13}$ with $0\le k \le n-1$. For pairs starting with odd indices, this means for each $k$ there is a triangle reduced sequence involving at least $R^{\circ}_{23}$ and having endpoints $x_{2k+1},x_{2k+2}$. If this is a proper type II reduced sequence, then it has to be a concatenation of loops. In particular $x_{2k+1}=x_{2k+2}$ or $(x_{2k+1},x_{2k+2})\in \Delta \subset R_0\subset R_1$, contradiction ($\Delta=\{(x,x)|x\in X\}$). Hence the triangle reduced sequence with endpoints $x_{2k+1},x_{2k+2}$ involves $R^{\circ}_{23}$ and at most {\bf one} of $R^{\circ}_{12},R^{\circ}_{13}$. Using absorption (prop. \ref{absorption}) we can thus build a proper triangle reduced loop of type I, involving $E^{\circ}_1,R^{\circ}_{23}$. This contradicts the definition of triangle join, so claim (\ref{l1}) holds.
\par Let $F_k=E_k\ast_{R_k}R$, and $F_k^{\circ}=F_k\setminus R$. Consider in $X$ a sequence
\[ x=x_0\stackrel{\tiny F_{i_0}}{\sim}x_1\stackrel{\tiny F_{i_1}}{\sim}x_2...\stackrel{\tiny F_{i_{n-2}}}{\sim}x_{n-1}\stackrel{\tiny F_{i_{n-1}}}{\sim}x_{n}=x,n\ge2.\]
We have to show that
\begin{equation}\label{l2}
\exists\quad k:\quad 1\le k<n\mbox{ so that } x_k\stackrel{\tiny R}{\sim}x_{k+1}.
\end{equation}
Suppose no such $k$ exists. Then for all $k$, $(x_k,x_{k+1})\in F_{i_k}^{\circ}=E_{i_k}\vee R\setminus R$. This means there exists a triangle reduced sequence $x_k=x_{k,1},x_{k,2},...,x_{k,n_k}=x_{k+1}$ with consecutive pairs alternating $E_{i_k}$ and $R$. Moreover, for each $k$, at least one $E_{i_k}^0$ must appear. Since the length of the initial sequence is $n\ge2$, at least two distinct $E^{\circ}_i$'s are involved. Using again absorption (prop. \ref{absorption}), we end up with a proper reduced triangle loop of type I, contradiction. So (\ref{l2}) has to be true and this ends the proof.
\rule{5 pt}{5 pt} 

\vspace{0.3 cm}
We have seen one instance when a triangle join that is not minimal is  actually a free amalgamated join. The next theorem gives a necessary and sufficient condition for a minimal triangle join to be a free join.
\begin{theorem}\label{freejoin}
Let $R=R_1\vee R_2\vee R_3=R_{12}\vee R_{23}\vee R_{13}$  be a minimal triangle join. Then
\[ R(=R_1\vee R_2)=R_1\ast_{R_{12}}R_2\Longleftrightarrow R_3=R_{13}\ast_{R_0}R_{23}. \]
\end{theorem}

\vspace{0.3 cm}
{\it Proof: } The fact that $R=R_1\vee R_2$ is clear from the definitions of $R,R_1,R_2$. \\
"$\Rightarrow:$" For $n\ge 1$ consider in $X$ the sequence
\[ x=x_0\stackrel{\tiny R_{13}}{\sim}x_1\stackrel{\tiny R_{23}}{\sim}x_2...\stackrel{\tiny R_{13}}{\sim}x_{2n-1}\stackrel{\tiny R_{23}}{\sim}x_{2n}=x.\]
Since $R_{13}\subset R_1$ and $R_{23}\subset R_2$, we also have
\[ x=x_0\stackrel{\tiny R_{1}}{\sim}x_1\stackrel{\tiny R_{2}}{\sim}x_2...\stackrel{\tiny R_{1}}{\sim}x_{2n-1}\stackrel{\tiny R_{2}}{\sim}x_{2n}=x.\]
Assume $(x_i,x_{i+1})\not\in R_0=R_{12}\cap R_{23}=R_{12}\cap R_{13}$ for all $0\le i<2n$.  Then $(x_i,x_{i+1})\not\in R_{12}$ for all $0\le i<2n$. This contradicts the fac that $R_1\vee R_2=R_1\ast_{R_{12}}R_2$.\\
"$\Leftarrow:$" For $n\ge 2$ consider in $X$ the sequence
\[ x=x_0\stackrel{\tiny R_{1}}{\sim}x_1\stackrel{\tiny R_{2}}{\sim}x_2...\stackrel{\tiny R_{1}}{\sim}x_{2n-1}\stackrel{\tiny R_{2}}{\sim}x_{2n}=x.\]
If $n=1$ it immediately follows that $(x_0,x_1)\in R_0$. Let $n\ge2$ and 
assume $(x_i,x_{i+1})\not\in R_{12}$ for all $0\le i<2n$. There exists then for each $i$ a triangle reduce sequence $x_i=y_{i,1},y_{i,2},...,y_{i,n_i}=x_{i+1}$ such that $n_i\ge 1$ and at least for one $k\in \{1,...,n_i\}$ we have
$x_{i,k}\stackrel{\tiny R^{\circ}_{13}}{\sim}x_{i,k+1}$ if $i$ is even, $x_{i,k}\stackrel{\tiny R^{\circ}_{23}}{\sim}x_{i,k+1}$ if $i$ is odd. Use absorption to end up with a triangle reduced loop alternating at least $R^{\circ}_{13}$ and $R^{\circ}_{23}$. If $R^{\circ}_{12}$ does not appear at all, i.e. $n_i=1$ for all $i$, then we contradict the hypothesis about $R_{13}$ and $R_{23}$. If $R^{\circ}_{12}$ is also used, we are dealing with a proper triangle reduced loop of type II, which, by the definition of the triangle join, has to be a concatenation of $R_i$-loops. The next claim is that
\begin{equation}\label{l3}
\mbox{ there is a triangle reduced sub-loop involving only }R^{\circ}_{13},R^{\circ}_{23}.
\end{equation}
Suppose on the contrary, that all the loops involve $R^{\circ}_{12}$. Note the following:
\begin{itemize}
\item because the loops  involve $R^{\circ}_{13} ,R^{\circ}_{12}$ or $R^{\circ}_{23},R^{\circ}_{12}$, their lengths, $n_i$, must be $n_i\ge4$. 
\item because the loops  can involve only $R^{\circ}_{13} ,R^{\circ}_{12}$ or $R^{\circ}_{23},R^{\circ}_{12}$, and we have assumed that there are no $R^{\circ}_{13},R^{\circ}_{23}$-loops, these loops will start and end around the initial $x_i$'s. Note that no loop could have endpoints some consecutive $x_i,x_{i+1}$ (or else $x_i=x_{i+1}$, contradiction). 
\item around a generic $x_i$ the following situations are possible: 
	\begin{itemize}
	\item[(i)] $x_{i-1,n_{i-1}-1}\stackrel{\tiny R^{\circ}_{23}}{\sim}x_i\stackrel{\tiny R^{\circ}_{12}}{\sim}x_{i,2}$\\
	or $x_{i-1,n_{i-1}-1}\stackrel{\tiny R^{\circ}_{12}}{\sim}x_i\stackrel{\tiny R^{\circ}_{23}}{\sim}x_{i,2}$
	\item[(ii)] $x_{i-1,n_{i-1}-1}\stackrel{\tiny R^{\circ}_{13}}{\sim}x_i\stackrel{\tiny R^{\circ}_{12}}{\sim}x_{i,2}$\\
	or $x_{i-1,n_{i-1}-1}\stackrel{\tiny R^{\circ}_{12}}{\sim}x_i\stackrel{\tiny R^{\circ}_{13}}{\sim}x_{i,2}$
	\item[(iii)] $x_{i-1,n_{i-1}-1}\stackrel{\tiny R^{\circ}_{12}}{\sim}x_i\stackrel{\tiny R^{\circ}_{12}}{\sim}x_{i,2}$
	\item[(iv)] $x_{i-1,n_{i-1}-1}\stackrel{\tiny R^{\circ}_{13}}{\sim}x_i\stackrel{\tiny R^{\circ}_{23}}{\sim}x_{i,2}$
	\end{itemize}
\item in cases (i),(ii), the loop ends/starts either at $x_{i-1,n_{i-1}-1}$ or $x_{i,2}$. For instance, if the loop ends at $x_{i,2}$, this forces $x_i\stackrel{\tiny R^{\circ}_{12}}{\sim}x_{i,2}$. In the other two cases the loop ends/starts at $x_i$.
 \end{itemize}

\par Let us look closer to what happens with the first loop. By the above remarks, the loop cannot end at 
$x_1$. Suppose it ends at $x_{0,n_0-1}$. Then $x_{0,n_0-1}=x_0=x$ and we must have $x_{0,n_0-1}\stackrel{\tiny R^{\circ}_{12}}{\sim}x_1$, hence $x_0\stackrel{\tiny R^{\circ}_{12}}{\sim}x_1$, contradiction. The other possibility again leads to the same conclusion. Hence $n_0=1$ and $(x_0,x_1)\in R^{\circ}_{13}$ which contradicts our assumption on the loops; (\ref{l3}) follows. In turn, this contradicts the fact that $R_{13}\vee R_{23}=R_{13}\ast_{R_0}R_{23}$. Hence the assumption that $(x_i,x_{i+1})\not\in R_{12}$ is false and this finishes the proof.
\rule{5 pt}{5 pt}

\section{\label{cost}Cost}

\subsection{\label{cost of groups}Cost of a graphing, SP1 equivalence relation and group}
\par In what follows we briefly review the notions of graphing, cost and results that we will use later on. For more details the reader is referred to \cite{GAB},\cite{KCH}. 
\par
A {\bf graphing} on $(X,\mu)$ is a countable family $\phi=\{\varphi_j:A_j\to B_j\}_{j\in J}$ of partial Borel isomorphisms between Borel subsets of $X$. Each $\varphi_j$ is called a {\bf generator}, with {\bf domain} $A_j$ and {\bf range} $B_j$.
A graphing $\phi$ generates an equivalence relation $R_{\phi}\subset X\times X$. It is the smallest equivalence relation satisfying $x\stackrel{\tiny R_{\phi}}{\sim}y$ if there exists $j\in J$ such that $x\in A_j$ and $y=\varphi_j(x)$.
Compositions of generators of $\phi$ are called $phi$-words. The graphing $\phi$ gives to each orbit $R_{\phi}[x]$ a natural structure of connected graph - the Cayley graph. If this graph is a tree for a.e. $x\in X$ we call $\phi$ a treeing.
\par  We say that the relations $R_1$ and $R_2$ on $(X_1,\mathcal{B}_1, \mu_1)$ and  $(X_2,\mathcal{B}_2, \mu_2)$ are isomorphic and denote it by $R_1\simeq R_2$, if there exists a Borel isomorphism $f:X_1\to X_2$ so that $f_{\ast}(\mu_1)=\mu_2$ and for $\mu_1$-a.e. $x\in X_1:\quad f(R_1[x])=R_2[f(x)]$. Here $f_{\ast}(\mu_1)$ is the measure induced by $f$ on $X_2$.
We say that $\phi$ is a graphing for the equivalence relation $R$ if $R_{\phi}\simeq R$.
\par A SP1 relation $R$ is called 
\begin{itemize}
\item {\bf finite} if for a.e. $x\in X$ the orbit $R[x]$ is finite. 
\item {\bf hyperfinite} if it is isomorphic to a non-decreasing limit of finite relations $R_n$: for a.e $x\in X, R_n[x]\subset R_{n+1}[x]$ and $\cup_{n\in \mathbb{N}}R_n[x]=R[x]$.
\end{itemize}
\begin{definition}Let $(X, \mathcal{ B},\mu)$ be a standard, non-atomic Borel probability space.
\begin{item}
\item[(1)] The cost of a graphing $\phi=\{\varphi_j:A_j\to B_j\}_{j\in J}$ on $(X,\mu)$ is the quantity:
\[ \mathcal{C}_{\mu}(\phi)=\sum_{j\in J}\mu(A_j).\]
\item[(2)] The cost of an SP1 relation $R$ on $(X,\mu)$ is the greates lower bound of the cost of graphings of $R$:
\[ \mathcal{C}_{\mu}(R)=\mbox{inf} \{\mathcal{C}_{\mu}(\phi)|R_{\phi}=R\}.\]
\item[(3)] The cost of a countable group $\Gamma$ is the greatest lower bound of the costs of relations $R_{\alpha}$, where $\alpha$ is a free measure preserving action of $\Gamma$ on some standard, non-atomic Borel probability space $(Y,\nu)$:
\[ \mathcal{C}(\Gamma)=\mbox{inf}\{ \mathcal{C}_{\mu}(R_{\alpha})|\alpha \mbox{ free action of } \Gamma \mbox{ on some } (Y,\nu), \nu(Y)=1\}.\]
\item[(4)] A group is said to be of fixed price if all the relations $R_{\alpha}$ have the same cost.
\end{item}
\end{definition}
%
The next theorem is Gaboriau's formula for computations of the cost of a particular kind of free join.
\begin{theorem}{\bf (\cite{GAB}, theorem IV.15)} Let $R=R_1\ast_{R_3}R_2$ be a SP1 relation on $(X,\mu)$, with $R_1,R_2$ of finite costs and $R_3$ hyperfinite. Then:
\[ \mathcal{C}(R)=\mathcal{C}(R_1)+\mathcal{C}(R_2)-\mathcal{C}(R_3).\]
\end{theorem}
\begin{remark}\label{r1} From Gaboriau's proof of the above theorem, we note that in general, if $R_3$ is a hyperfinite SP1 subrelation of the finite cost SP1 relations $R_1,R_2$ on $(X,\mu)$, then
\[  \mathcal{C}(R_1\vee R_2)\le \mathcal{C}(R_1)+\mathcal{C}(R_2)-\mathcal{C}(R_3). \]
\end{remark}

There is also a converse to the previous theorem, due also to Gaboriau:
\begin{theorem}{\bf (\cite{GAB}, theorem IV.16)} Let $R_1,R_2$ be two SP1 relations of finite cost on $(X,\mu)$, let $R=R_1\vee R_2$ be the SP1 relation  generated by $R_1,R_2$ and let $R_3$ be a finite subrelation of $R_1,R_2$. \\ If $\mathcal{C}(R)= \mathcal{C}(R_1)
+ \mathcal{C}(R_2)- \mathcal{C}(R_3)$, then $R=R_1\ast_{R_3}R_2$.
\end{theorem}

\subsection{\label{cost 2}Cost of triangle amalgams} Using Gaboriau's results on the cost of a free join, we can now give an estimate of the cost of a particular kind of triangle join.
\begin{proposition}
Let $R$ be the minimal triangle join of the finite cost SP1 relations $R_1,R_2,R_3$ on $(X,\mu)$. Assume further that $R_i\cap R_j=:R_{ij}$ are finite SP1 relations and that none of the vertex relations is a free join of the edges. Then
\[\mathcal{C}(R)< min\{\mathcal{C}(R_i)+\mathcal{C}(R_j)-\mathcal{C}(R_{ij})|1\le i<j\le3\}.\]
\end{proposition}
{\it Proof: } Clearly $R=R_1\vee R_2=R_2\vee R_3=R_1\vee R_3$. By the above remark \ref{r1}, we get that 
\[ \mathcal{C}(R)\le min\{\mathcal{C}(R_i)+\mathcal{C}(R_j)-\mathcal{C}(R_{ij})|1\le i<j\le3\}. \]
Since none of the vertex relations is a free join of the edges, using th. \ref{freejoin}  we see that none of the  relations $R_1\vee R_2=R_2\vee R_3=R_1\vee R_3$ is a free join, hence the inequalities are strict and the result follows.
\begin{corollary}\label{cost}
Suppose $G$ is the minimal triangle join of the amenable (finite or countable) $G_{12},G_{13},G_{23}$. Assume further that $G_1,G_2$ and $G_3$ are of fixed finite cost. Then
\[ \mathcal{C}(G)<min\{ \mathcal{C}(G_i)+ \mathcal{C}(G_j)-\mathcal{C}(G_{ij})|1\le i<j\le3\}.\]
\end{corollary}
\subsection{An application:  the cost of  Thompson's group $G_{2,1}$}
We follow \cite{Brown} to give an estimation of the cost of Thompson's infinite simple group $G_{2,1}$. Let $\mathcal{G}$ be the family of finitely presented infinite simple groups introduced by Higman, generalizing R.J. Thompson's group of dyadic homeomorphisms of the Cantor set. In \cite{Brown} the main theorem states:
\begin{theorem} For any $G\in \mathcal{G}$ and any $n\ge 1$ there is an $(n-1)$-connected $n$-dimensional simplicial complex $X$ such that $G$ acts on $X$ with finite stabilizers and with an $n$-simplex as fundamental domain.
\end{theorem}
As a consequence, theorem 3 in \cite{Brown} states that, for $n=2$, Thompson's group $G_{2,1}$ has the following combinatorial presentation as a realizable triangle amalgam,
\begin{table}[ht]
\begin{center}
\begin{tabular}{ccccccccc}

$ $&$ $&$ $&$ $&$S_{p+2}$&$ $&$ $&$ $&$ $\\
$ $&$ $&$ $& $\nearrow $ &$ $&$\nwarrow $&$ $&$ $&$ $\\
$ $&$ $&$S_{p-2}\times S_2 $&$ $&$ $&$ $&$S_p $&$ $&$ $\\
$ $&$ $&$ $&$ $&$ $&$ $&$ $&$ $&$ $\\
$ $&$\swarrow $&$ $&$ $&$$&$ $&$ $&$\searrow $&$ $\\
$ $&$ $&$ $&$ $&$ $&$ $&$ $&$ $&$ $\\
$S_p $&$ $&$\longleftarrow$&$ $&$S_{p-1}$&$ $&$\longrightarrow$&$ $&$S_{p+1} $\\
\end{tabular}
\end{center}
\end{table}
where $S_r$ stands for the group of permutations of ${1,2,...,r}$. The inclusions, as constructed by Brown, are (bellow we present the general case, which was only mentioned to work; Brown explicitly constructed the inclusions only for $p=5$, but he mentions that one gets similar triangles for every $p\ge 5$):
\begin{itemize}
\item $S_{p-1}\hookrightarrow S_r$ for $r=p,p+1$ are the standard inclusions that permute $\{1,2,...,p-1\}\subset \{1,...,r\}$;
\item $S_p\hookrightarrow S_r$ for $r=p,p+1$ are the inclusions that permute $\{1,2,3,r-1,r\}\subset \{1,...,r\}$;
\item $S_{p-2}\times S_{2}$ is embedded in $S_p$ so that the first factor permutes $\{1,2,...,p-2\}$ and the second factor acts on $\{p-1,p\}$;
\item $S_{p-2}\times S_2\hookrightarrow S_{p+2}$ so that the first factor permutes $\{1,2,...,p-2\}$ and the second factor sends the nontrivial permutation not into a transposition, but rather to the product $(p-1,p+1)(p,p+2)$ of two transpositions.
\end{itemize}
Note that for every $p\ge 5$, the triangle amalgam yields the same group $G_{2,1}$.
We would like to emphasize that this realizable amalgam is a {\bf minimal} one, i.e. the vertex groups are generated by the images of the edge groups. Since none of the vertex groups is a free product, the induced relations on any $(X,\mu)$ cannot be free joins. From Gaboriau's price list in \cite{GAB}, we know that a finite group $H$ has cost $\mathcal{C}(H)=1-\frac{1}{|H|}$. Hence the above corollary \ref{cost} gives the following estimate:
\[\mathcal{C}(G_{2,1})<min\{1-\frac{1}{p!}-\frac{1}{(p+2)!}+\frac{1}{2!\cdot (p-2)!},1-\frac{1}{p!}-\frac{1}{(p+1)!}+\frac{1}{(p-1)!},\] \[1-\frac{1}{(p+1)!}-\frac{1}{(p+2)!}+\frac{1}{(p-1)!}\}=\]
\[=1-\frac{1}{(p+1)!}-\frac{1}{(p+2)!}+\frac{1}{p!},\mbox{ for all }p\ge 5.\]
Since $G_{2,1}$ is infinite, it must have cost greater or equal to $1$ (\cite{GAB}, properties VI.3.). Letting $p\to \infty$ we obtain:
\[ \mathcal{C}(G_{2,1}) =1.\qquad \qquad \rule{5pt}{5pt}\]
\begin{remark} The groups $F_{2,1},G_{2,1}$ (also another one, called $T_{2,1}$) were  originally introduced by Thomspon in 1965 and provided the first examples of finitely presented infinite simple groups. It is known that $G_{2,1}$ contains, for every positive $n$, the symmetric group on $n$ letters. As a consequence, every finite group can be embedded in $G_{2,1}$ (\cite{Brown}). Moreover, every countable locally finite group can be embedded in $G_{2,1}$ (\cite{GH}).

One open question about these groups still stands: are Thompson's groups $F_{2,1}$,$G_{2,1}$ amenable? Of course, amenability of $G_{2,1}$ would imply amenability for $F_{2,1}$. Gaboriau's price list shows that infinite amenable groups all have cost $1$ and are treeable. 
\end{remark}


\begin{acknowledgements} I want to thank my advisor, Florin R\u adulescu,  for the constant help and guidance. 
\end{acknowledgements}

 \end{document}